\documentstyle[epsf,12pt]{article}
\normalbaselines
\catcode `\@=11

\def\section{\@startsection {section}{1}{\z@}{-3.5ex plus -1ex minus
     -.2ex}{2.3ex plus .2ex}{\normalsize\bf}}
\def\subsection{\@startsection{subsection}{2}{\z@}{-3.25ex plus -1ex minus
 -.2ex}{1.5ex plus .2ex}{\normalsize\bf}}

\def\thebibliography#1{\section*{References\markboth
  {REFERENCES}{REFERENCES}}\list
  {[\arabic{enumi}]}{\settowidth\labelwidth{[#1]}\leftmargin\labelwidth
  \advance\leftmargin\labelsep
  \usecounter{enumi}}
  \def\newblock{\hskip .11em plus .33em minus -.07em}
  \sloppy
  \sfcode`\.=1000\relax}
 
\catcode `\@=12

\newtheorem{theorem}{Theorem}

\begin{document}

\begin{flushright} Preprint DM/IST 10/97
\end{flushright}
\vspace*{1cm}
\begin{center} { \bf KNOT INVARIANTS FROM A KOHNO-KONTSEVICH INTEGRAL FOR
BRAIDS}
\vspace{1cm}\\  Roger Picken$^{1}$\\ \vspace{0.3cm}
${}^{1}$  Departamento de Matem\'{a}tica \& Centro de Matem\'{a}tica 
Aplicada,\\ Instituto Superior T\'{e}cnico, Av. Rovisco Pais,\\ 1096 Lisboa
Codex, Portugal.
\end{center}

\vspace*{0.5cm}

\begin{abstract}
\noindent 
An invariant of knots is constructed from an integral for geometric braids
due to Kohno and Kontsevich. It takes values in a quotient by a certain
ideal of the algebra generated by chord diagrams over the circle. 

\end{abstract}

\section{\hspace{-4mm}.\hspace{2mm} Preliminaries}

\hspace*{0.8cm}
In this introductory section we present some basic notions concerning
braids and knots, as well as the definition and a couple of properties of
an integral for braids due to Kohno and Kontsevich. For more detailed
reviews of these subjects, containing further background material, the
reader is referred to articles by Birman \cite{BirAMS} and Bar-Natan
\cite{BN}.

Let $0<x_1<\cdots <x_n$ be a set of $n$ numbers lying on the positive real
axis of the  complex plane $\bf C$. Identify ${\bf R}^3$ with ${\bf
C}\times{\bf R}$,  coordinatized by $(z,t)$ with $z\in{\bf C}$ and $t\in
{\bf R}$. An $n$-strand braid, or  simply braid, $b$, is a one-dimensional
piecewise-smooth oriented submanifold of 
${\bf R}^3$ lying between the planes $t=0$ and $t=1$, which is topologically
the  disjoint union of $n$ intervals called the strands of $b$, whose
boundary is 
$\left\{(x_i,j):i=1,\dots,n,\,j=0,1\right\}$, and whose strands point
upwards everywhere  in the sense that the $t$ component of the orientation
vector is everywhere positive.  Frequently one makes the ``natural'' choice
of endpoints $x_i=i$ for $i=1,\ldots,n$,  but for reasons which will become
clear we wish to allow greater flexibility for the  endpoints.

Two $n$-strand braids $b_1$ and $b_2$ are said to be equivalent if $b_1$ is
carried  into $b_2$ under a diffeomorphism of ${\bf C}\times  [0,1]$ which
preserves each horizontal plane, as well as the real axes of the top  and
bottom planes, and is identity-connected when restricted to the real axes
of the top  and bottom planes. In the special case of the diffeomorphism
restricting to the identity  on the top and bottom planes, we will speak of
restricted equivalence of braids. The  unrestricted equivalence allows the
endpoints to move within the positive real axes  without changing their
order.

Two $n$-strand braids $b_1$ and $b_2$ sharing the same endpoints may be 
multiplied: $b_1b_2$ is the braid obtained by shrinking both braids by a
factor $1/2$ in  the $t$ direction and putting $b_2$ on top of $b_1$. Given
a braid $b$, its inverse 
$b^{-1}$ is the mirror reflection of $b$ in the bottom plane, lifted in the
positive $t$  direction to lie again between the $t=0$ and $t=1$ planes. The
restricted equivalence  classes of braids sharing a fixed set of endpoints
form a group under these operations,  with the unit being the class of the
trivial braid whose strands point vertically upwards  everywhere. 

Braids may be used as a tool to study knots and links through the closure
construction.  Since $b$ is compact it is contained in $U\times  [0,1]$ for
some open neighbourhood $U\subset {\bf C}$. The closure of
$b$ is  the knot or link obtained by connecting the corresponding points
$(x_i,0)$ and 
$(x_i,1)$ for $i=1,\ldots,n$ by $n$ parallel curves contained in the
complement of 
$U\times ]0,1[$. Since knots and links are identified up to ambient isotopy
the closure is  well-defined. Any knot or link can be obtained as the
closure of some braid by a  theorem of Alexander \cite{Al}. Furthermore the
closure of two equivalent braids gives  the same knot or link. When two
braid equivalence classes have the same closure they  are called Markov
equivalent. Markov's theorem \cite{Ma}\cite{Bi} states that Markov 
equivalence is generated by two types of move: (Markov 1) $b_1b_2
\rightleftharpoons  b_2b_1$ and (Markov 2) $b \rightleftharpoons (b
\amalg |){\sigma_n}^{\pm 1}$, where $b$  is an $n$-strand braid contained
in $U\times  [0,1]$ for some open neighbourhood $U\subset {\bf C}$, $ b
\amalg |$ is  an $(n+1)$-strand braid consisting of $b$ together with an
$(n+1)$th strand of the  form $\left\{(x_{n+1},t):\,0\leq t\leq 1\right\}$
contained in the complement of 
$U\times  [0,1]$, and ${\sigma_n}^{\pm 1}$ is the elementary $(n+1)$-strand
braid where the 
$n$th and $(n+1)$th strands cross over ($+$) or under ($-$) each other,
whilst the other  strands are vertical. Two strands crossing over (under) is
taken to mean that they circle  halfway round each other in an anticlockwise
(clockwise) direction as $t$ goes from $0$ to 
$1$.

\begin{figure}[h]
\begin{center}
\mbox{\epsfxsize=1\hsize\epsfbox{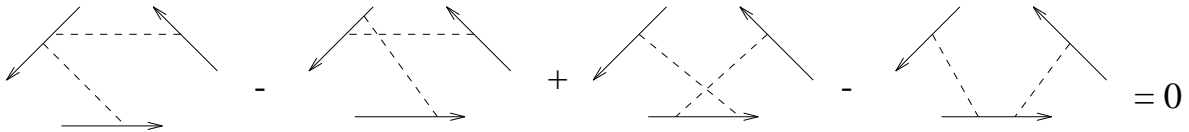}}
\caption{The 4T relations}
\label{fig:4T}
\end{center}
\end{figure}

The notion of chord diagrams for one-dimensional manifolds has been used by
a number of  authors, in particular in the context of Vassiliev invariants 
\cite{Vas1}
\cite{Vas2}
\cite{Kon}
\cite{BN}
\cite{Le-Mu1}
\cite{Ka-Tu}. 
Let $X$ be an
oriented one-dimensional manifold.  A chord diagram over $X$ is a finite
set of disjoint pairs of distinct points belonging to the  interior of $X$,
represented pictorially by dashed lines connecting the points of each pair. 
Chord diagrams are identified up to orientation and component-preserving
diffeomorphisms  of
$X$. The space of chord diagrams over $X$, denoted ${\cal A}(X)$, is the
vector space  generated over $\bf C$ by the chord diagrams over $X$,
quotiented by the four-term (4T)  relations. The 4T relations hold between
any four diagrams which are identical apart from in  three disjoint open
intervals contained in
$X$, where they are as in Figure~\ref{fig:4T}.

 The vector space ${\cal A}(X)$ is naturally graded over $\bf N$ by the
number of  chords. Denote by ${\cal A}_m(X)$ the subspace of ${\cal A}(X)$
generated by chord  diagrams with precisely $m$ chords, modulo the 4T
relations. We will in fact be  considering the completion of ${\cal A}(X)$
by this grading, also denoted ${\cal  A}(X)$. For the case $X=S^1$ we write
$\cal A$ instead of ${\cal A}(S^1)$. Two chord  diagrams over $S^1$ can be
multiplied by forming the connected sum of the two  circles, making the
connection in regions free of chord endpoints. Because of the 4T  relations
this product is well-defined (independent of how the connected sum is 
formed) and, extending the product to the whole of $\cal A$ by linearity, 
$\cal A$  becomes a commutative algebra with unit (namely the $0$-chord
diagram over $S^1$  denoted 
$\bigcirc$). In the case where $b_1$ and $b_2$ are $n$-strand braids there
is a natural product  operation $ {\cal A}(b_1)\times{\cal
A}(b_2)\rightarrow{\cal A}(b_1b_2)$ given, for  two individual chord
diagrams, by putting the chord diagram over $b_2$ on top of the  chord
diagram over $b_1$, and extended by linearity to general elements.  Finally
for two manifolds $X_1$, $X_2$, there is a product ${\cal A}(X_1)\times
{\cal  A}(X_2)\rightarrow {\cal A}(X_1
\amalg X_2)$ $(c_1,c_2)\mapsto c_1\otimes c_2$, where
$\amalg$ is the disjoint union, induced by the natural inclusion maps from
the manifolds to their  disjoint union.

The  following integral associated with a braid $b$ and taking values in
${\cal A}(b)$ was  introduced by Kontsevich 
\cite{Kon} in the context of Vassiliev invariants and developed further in 
\cite{BN}
\cite{Le-Mu2}. It is closely related to an integral introduced
earlier by Kohno \cite{Koh1}
\cite{Koh2}
 taking values in a matrix algebra and reduces to the Kohno integral in a 
matrix representation of ${\cal A}(b)$. Take $m$ values of
$t$, $0<t_1<\cdots  <t_m<1$ and for each $t_i$ choose a pair of distinct
points $(z_i(t_i),t_i)$ and 
$(z^\prime_i(t_i),t_i)$ belonging to $b$. This is called a pairing, $P$, and
there is an  obvious chord diagram $D_P\in {\cal A}(b)$ associated with it.
Now we define:
\begin{equation} \label{eq:KI} Z(b)=\sum_{m=0}^{\infty}\frac{1}{(2\pi
i)^m}\int_{0<t_1<\cdots <t_m<1}
\sum_{{\rm pairings}\,P}\prod_{i=1}^{m}\frac{d(z_i-z^\prime_i)}{
z_i-z^\prime_i } D_P.
\end{equation}
$Z(b)$ is best understood geometrically as a path-ordered exponential
representing  parallel transport in the configuration space of $n$ particles
in the plane ${\bf C}^n
\setminus \Delta$ where $\Delta=\left\{(z_1,\ldots, z_n):z_i=z_j\, {\rm for}
\,
i\neq  j\right\}$, with respect to the flat Knizhnik-Zamolodchikov
\cite{KZ} connection $\omega_{\rm KZ}=(1/2\pi
i)\sum_{i<j=1}^{n}d(z_i-z_j)/(z_i-z_j)  H_{ij}$ where $H_{ij}\in {\cal
A}(b)$ is the chord diagram with a single chord  connecting strands $i$ and
$j$. Due to the flatness of $\omega_{\rm KZ}$, $Z(b)$ has  the property:
\begin{theorem} $Z(b)$ respects the restricted equivalence of braids.
\end{theorem} Furthermore   it has the
multiplicative property \cite{BN}\cite{Le-Mu1}: 
\begin{theorem} $Z(b_1b_2)=Z(b_1)Z(b_2)$.
\end{theorem}

\section{\hspace{-4mm}.\hspace{2mm}Construction of a knot invariant}

\hspace*{0.8cm} Let us for simplicity restrict attention from now on to
those braids whose closure is a  knot, i.e. a one-component link. The
construction which follows may be easily adapted  to the general case. For
such braids $b$ there is a natural map $p:{\cal  A}(b)\rightarrow \cal A$
given by the inclusion of $b$ in its closure, topologically a  circle. $pZ$
enjoys the following two properties:
\begin{theorem}$pZ(b_1b_2)=pZ(b_2b_1)$.
\end{theorem} 
{\em Proof:} This follows from Thm.2 and the obvious property 
$p(c_1c_2)=p(c_2c_1)$ for $c_1\in {\cal A}(b_1)$, $c_2\in {\cal A}(b_2)$.
$\Box$
\begin{theorem} $pZ$ respects the equivalence of braids.
\end{theorem}  
{\em Proof:} The Kontsevich integral is defined also for braid-like
objects having 
$n$ upward-pointing strands but whose endpoints in the top and bottom
planes, both  identified with $\bf C$, are not the same. Let $b_1$ and $b_2$
be equivalent braids  such that the endpoints of $b_1$ and $b_2$ are 
$\left\{(x_i,k):i=1,\ldots,n,\,k=0,1\right\}$ and 
$\left\{(x^\prime_i,k):i=1,\ldots,n,\,k=0,1\right\}$ respectively. Consider
the braid-like  object $\tilde{b}$ which has bottom endpoints
$\left\{(x_i,0):i=1,\ldots,n\right\}$ and top endpoints
$\left\{(x^\prime_i,1):i=1,\ldots,n\right\}$, and whose strands are 
straight lines connecting corresponding endpoints. Let $\tilde{b}^{-1}$ be
its inverse in  the sense of the mirror construction used for braids. The
equivalence of $b_1$ and 
$b_2$ implies that $b_1$ and $\tilde{b}b_2\tilde{b}^{-1}$ are equivalent in
the restricted  sense. Adapting Thm.~3 to braid-like objects one has
$pZ(b_1)=pZ(\tilde{b}b_2\tilde{b}^{-
1})=pZ(b_2\tilde{b}^{-1}\tilde{b})=pZ(b_2)$, where the last equality uses
Thm.~1. $\Box$

Now we calculate $pZ$ for the simplest non-trivial $2$-strand braid
$\sigma_1$ and its  inverse:
\begin{equation} 
pZ({\sigma_1}^{\pm 1})= \bigcirc + \sum_{m=1}^{\infty}(\pm 1/2)^m/m!\,
\bigcirc _m
\end{equation} where $\bigcirc _m$, for $m=1,\ldots,\infty$, denotes the
$m$-chord diagram over
$S^1$ such that the  endpoints of each chord are diametrically opposite
points of the circle. This is because  the coefficient of  $\bigcirc _m$
coming from equation~(\ref{eq:KI}) reduces to a path-ordered integral along
a single contour $C$ given by $z(t)=\exp (\pi it)$ for $0\leq t\leq 1$,
namely
$(1/2\pi i)^m {\cal P}\int_C dz_1/z_1\cdots dz_m/z_m$  where $\cal P$
denotes path-ordering, which is given by:
\begin{equation}
(1/2)^m\int_{t_m=0}^1\int_{t_{m-1}=0}^{t_m}\ldots\int_{t_1=0}^{t_2}
dt_1dt_2\ldots dt_m = (1/2)^m 1/m!\,.
\end{equation} Now we define the $\cal A$ elements
\begin{equation} r_\pm = pZ({\sigma_1}^\pm) - 
\bigcirc .
\end{equation}
$r_\pm$ are non-invertible and thus the ideal ${\cal I}=r_+{\cal A}+r_-{\cal
A}$ is  non-trivial. Let $\tilde{\cal A}$ denote the quotient of $\cal A$ by
this ideal and let 
$k:{\cal A}\rightarrow \tilde{\cal A}$ be the canonical projection.

\begin{figure}[tbp]
\begin{center}
\mbox{\epsfxsize=0.3\hsize\epsfbox{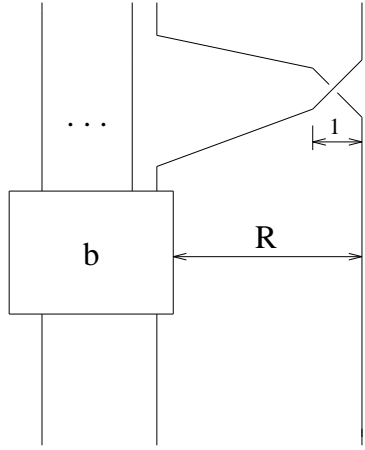}}
\caption{Braid used in the proof of Thm.~5}
\label{fig:bR}
\end{center}
\end{figure}

The following main theorem employs a limit technique similar to one used in
Le and Murakami
\cite{Le-Mu1}, Lemma 2.3.4. 
\begin{theorem}\label{th:main}
 Let a knot $K$ be presented as the closure of a braid $b$. Then 
$Y(K)=kpZ(b)$ is well-defined and thus gives rise to an $\tilde{\cal
A}$-valued knot  invariant.
\end{theorem} 
{\em Proof:} Since $pZ$ respects braid equivalence and the first
Markov move (Thms.~3 and  4) it remains to show that $kpZ$ respects the
second Markov move. Let $b$ be an 
$n$-strand braid and $(b
\amalg |){\sigma_n}^{\pm 1}$ be the corresponding $(n+1)$-strand braid in
the move.  Translate the final strand and the over/undercrossing
${\sigma_n}^{\pm 1}$ out to the  right by a distance $R-1$ as in
Figure~\ref{fig:bR}.
 Call this new braid $b_n(R)$. It factorizes as
\begin{equation} 
b_n(R)=(b
\amalg_R |) \tilde{b}_R 
(|^{n-1}\amalg_R {\sigma_1}^{\pm 1}) 
\tilde{b}^{-1}_R
\end{equation} 
where $b \amalg_R |$ is the braid $b \amalg |$ with a distance $R$
separating $b$ and the final strand, ${\tilde{b}_R}$ is a
braid-like object with all strands pointing vertically upwards except the
$n$th one, $|^{n-1}\amalg_R {\sigma_1}^{\pm 1}$ is the elementary
braid ${\sigma_n}^{\pm 1}$ where the final two strands have been moved out
to the right by a distance $R$ 
and finally ${\tilde{b}_R}^{ -1}$ is the inverse of
${\tilde{b}_R}$
(see  Figure~\ref{fig:bR}).
 In the corresponding factorization of $Z(b_n(R))$ we split each factor into
an $R$-independent and an $R$-dependent part:
\begin{equation} \label{eq:fact1} 
Z((b
\amalg_R |)=Z(b)\otimes Z_0(|) +O(1/R)
\end{equation}
\begin{equation} \label{eq:fact2} 
Z(|^{n-1}\amalg_R {\sigma_1}^{\pm 1})=Z_0(|^{n-1})\otimes
Z({\sigma_1}^{\pm 1})+O(1/R)
\end{equation}
\begin{equation} \label{eq:fact3} 
Z(\tilde{b}^{\pm 1}_R)
=Z_0(\tilde{b}^{\pm 1}_R) + \sum_{i\geq 1}
Z_i(\tilde{b}^{\pm 1}_R)=0
\end{equation} 
The second terms on the right hand sides of 
(\ref{eq:fact1}) and
(\ref{eq:fact2}) are $O(1/R)$ as $R\rightarrow\infty$ since the
corresponding chords are multiplied by  coefficient integrals of order
$O(\ln (R+1) -\ln R)\sim O(1/R)$. In 
(\ref{eq:fact3}) the $R$-independent term is $ Z_0(\tilde{b}^{\pm 1}_R)$
(note that
${\cal  A}_0(\tilde{b}^{\pm 1}_R)$ is $R$-independent) since for $R=1$, $
\tilde{b}^{\pm 1}_R$  is the trivial braid and thus $\sum_{i\geq 1}
Z_i(\tilde{b}^{\pm 1}_R)=0$. The terms of chord number $i$ in this
remainder are
$O((\ln R)^i)$ as
$R\rightarrow
\infty$. Now $pZ(b_n(R))$ is $R$-independent by  Thm.~4. Since the only
convergent limits which can be obtained by multiplying terms of  order
$O(1/R)$ and terms of order $O((\ln R)^i)$ are $0$, and since $pZ(b_n(R))$
as a  whole is convergent, we conclude that $pZ((b
\amalg |){\sigma_n}^{\pm 1})$ is given by
\begin{equation} pZ((b
\amalg |){\sigma_n}^{\pm 1})= p((Z(b)\otimes Z_0(|))Z_0(\tilde{b}_R)  
(Z_0(|^{n-1} )\otimes Z({\sigma_1}^{\pm 1}))Z_0({\tilde{b}}^{-1}_R))
\end{equation} When $k$ is applied to both sides of this equation, on the
right-hand-side only the $0$- chord term remains in the third factor and we
thus obtain
$ Y((b
\amalg |){\sigma_n}^{\pm 1})=Y(b).
$ $\Box$

We finally prove a small non-triviality result. 
\begin{theorem} The invariant $Y$ distinguishes the trefoil knot and its
mirror image.
\end{theorem} 
{\em Proof:} The trefoil may be presented as the closure of
${\sigma_1}^3$, whereas its  inequivalent mirror image is the closure of
${\sigma_1}^{-3}$. Now we expand 
$Y({\sigma_1}^3)$ and $Y({\sigma_1}^{-3})$ up to 3-chords:
\begin{eqnarray*} Y({\sigma_1}^3) &=&
\bigcirc +(3/2)
\bigcirc _1 +(3/2)^2 1/2!
\bigcirc _2 + (3/2)^3 1/3!
\bigcirc _3 + \cdots\\ &=& 
\bigcirc + 6r_+ + 3r_- + (1/2)
\bigcirc _3 + \cdots
\end{eqnarray*} and by a similar calculation
\begin{eqnarray} Y({\sigma_1}^{-3}) &=&
\bigcirc + 3r_+ + 6r_- - (1/2)
\bigcirc _3 + \cdots
\end{eqnarray} where the dots represent terms with higher chord diagram
number. An arbitrary linear combination of the form 
$\bigcirc _3+\cdots$ does not belong to the ideal $\cal I$, as may be
easily established by  studying 4T relations for $3$-chord diagrams in $\cal
A$. Thus the coefficient of $\bigcirc _3$, being the first non-trivial
coefficient in the $Y$ invariant, differs for the trefoil  and its mirror
image. $\Box$

We conclude with a few remarks. 

Previous approaches to using the Kontsevich integral for obtaining knot and
link  invariants \cite{BN}\cite{Kon}\cite{Le-Mu1} led to severe
computational problems. The integral~(\ref{eq:KI}), adapted to tangles, has
to be multiplied by a certain element of
$\cal A$ before it  becomes a knot invariant. This element is rich in
structure, being related to Euler  sums/Zagier zeta functions 
\cite{Le-Mu1}\cite{BBB}\cite{Zag},
 but complicated to use. Whilst the integral~(\ref{eq:KI}) for braids
remains far from trivial, one may hope that the fairly rigid structure of
braids  will enable efficient calculational methods to be found.

The invariant takes values in $\tilde{\cal A}={\cal A}/{\cal I}$, a
satisfyingly large space.  The ideal $\cal I$ is, loosely speaking, only
twice as large as the ideal $
\bigcirc _1\cal A $ of separated chord diagrams, which occurs in the
Bar-Natan approach to  Vassiliev invariants
\cite{BN}. It should be pointed out that in our construction there is no
need to quotient out by this  ideal.

\section{\hspace{-4mm}.\hspace{2mm}Acknowledgements}

\hspace*{0.8cm}
I am grateful to Nenad Manojlovic for help with the figures.

\end{document}